\documentclass{article}

\usepackage{microtype}
\usepackage{graphicx}
\usepackage{float}
\usepackage{subfigure}
\usepackage{booktabs} 

\usepackage{mbubble}
\mbubblestyle{ms}{color=blue,font={\tiny\sf}}

\usepackage{amsmath}
\usepackage{amssymb}
\usepackage{varioref}
\usepackage{hyperref}
\usepackage{cleveref}
\usepackage{graphicx}
\usepackage[utf8]{inputenc}
\usepackage{caption}

\usepackage{tabularx}
\usepackage{array}
\usepackage{booktabs}
\usepackage{longtable}

\usepackage{hyperref}

\usepackage[accepted]{icml2018}

\newcommand{\tr}{\mathsf{T}} 
\newcommand{\E}{\mathbb{E}} 
\newcommand{\given}{\ |\ } 
\newcommand{\toto}{{t-1,t-1}}
\newcommand{\tot}{{t-1,t}}
\newcommand{\ttt}{{t,t}}
\newcommand{\xh}{\hat{x}}
\newcommand{\bK}{\bar{K}}

\icmltitlerunning{Empirical fixed point bifurcation analysis}

\begin{document}

\twocolumn[
\icmltitle{Empirical fixed point bifurcation analysis}



\icmlsetsymbol{equal}{*}

\begin{icmlauthorlist}
\icmlauthor{Gergo Bohner}{gatsby}
\icmlauthor{Maneesh Sahani}{gatsby}
\end{icmlauthorlist}

\icmlaffiliation{gatsby}{Gatsby Unit, UCL, London, United Kingdom}

\icmlcorrespondingauthor{Gergo Bohner}{gbohner@gatsby.ucl.ac.uk}

\icmlkeywords{Machine Learning, ICML}

\vskip 0.3in
]



\printAffiliationsAndNotice{}  

\begin{abstract}
In a common experimental setting, the behaviour of  a noisy dynamical system is monitored in response to manipulations of one or more control parameters. Here, we introduce a structured model to describe parametric changes in qualitative system behaviour via stochastic bifurcation analysis. In particular, we describe an extension of Gaussian Process models of transition maps, in which the learned map is directly parametrized by its fixed points and associated local linearisations. We show that the system recovers the behaviour of a well-studied one dimensional system from little data, then learn the behaviour of a more realistic two dimensional process of mutually inhibiting neural populations.
\end{abstract}

\section{Introduction}
\label{Introduction}







Data, especially from biological systems, is often gathered through slow, noisy measurements, and as such standard modelling tools are often not suited to processing the raw experimental data. Time evolution of processes is frequently ignored, both on short time scales, such as during a single measurement, and on long ones -- through repeated trials taking place over weeks or months. 

Here we wish to address two features of experiments that do not always make the cut. Firstly, we explicitly model the short term evolution of the system using a stochastic dynamical system. Secondly, we wish to understand, how the behaviour changes with respect to varying a parameter that the experimenter either has direct access to, or that changes slowly and may be measured over long time periods.

Unfortunately, providing simple summaries and comparisons of dynamical systems learned from data is no easy task, and is a very active research area \cite{Sussillo2013, Nonnenmacher2017, Sussillo2016}. One popular approach has been to use a simple model to fit to the data, which enabled scientists to glean insight into how systems change and evolve \cite{Park2015}. A frequent criticism of this approach is that they do not adequately capture the data itself, so one has to be very careful when interpreting results.

On the other hand, more complex models may provide very good description of almost arbitrary data, especially with the recent rise of various deep learning algorithms. Despite their success in predicting newly observed data, we are far from an understanding of how changes in fitted parameters of such systems relate to changes in the data.

We take a middle approach here: We follow changes of dynamical systems via the bifurcations of their fixed points as the experimental parameter varies, a type of analysis mostly applied to known, continuous time and noiseless systems. In order to use it in the wild, we need to relax these constraints, and use a non-parametric stochastic map model, in which the fixed points can be identified robustly. For this purpose, we employ a Gaussian Process model of the map, which we parametrise in terms of its fixed points and associated Jacobians.

The paper is organised as follows: In \cref{sec:bifurc} we review deterministic bifurcation analysis on a simple example, and discuss the implications of noise. In \cref{sec:model}, we first describe how to infer a map from time series data using Gaussian Processes, then add the fixed points and derivative structures to our model of the map. Finally, in \cref{sec:experiments}, we prove the feasibility of our approach on a one-dimensional pitchfork bifurcation, then examine a more realistic experimental scenario using a model of mutually inhibiting neural populations.

To provide a concise description of the model with the added clarity of indexing, we employ the indicial notation with Einstein summation convention throughout the paper. In short, the number of lower indices describes the order of a tensor, and repeated indices in the same term induces summation, e.g. $\mathbf{a}^\tr\mathbf{b} \triangleq a_ib_i = \sum_ia_ib_i$, and $\mathbf{A}\mathbf{b} \triangleq A_{ij}b_j$. Furthermore, the third order Kronecker delta tensor $\delta_{ijk}$ is $1$ iff. $i=j=k$. For a much more detailed description of the notation system, we refer the reader to Appendix A.

\begin{figure*}[t]
\centering
\includegraphics[width=0.7\linewidth]{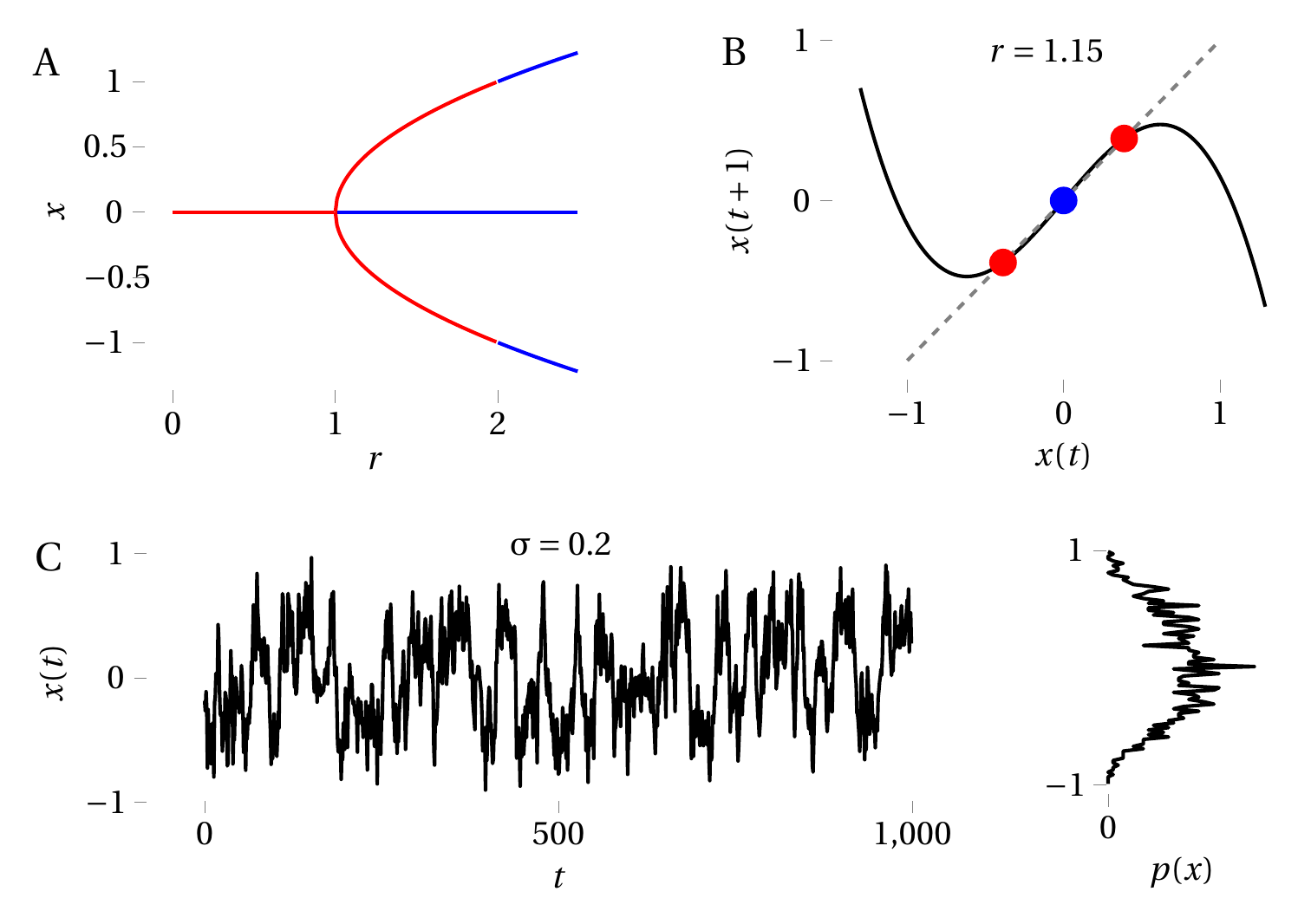}
\caption{Bifurcation analysis concepts. {\it A.} Pitchfork bifurcation diagram for \cref{eq:bifurc}. Throughout the paper red and blue colours indicate stable and unstable fixed points, respectively. {\it B.} The transition map visualised for $r=1.15$. The dashed line indicates potential fixed points $x(t)=x(t+1)$, and the circles indicate true fixed points. {\it C.} Realisation of a single noisy time series using $r=1.15, \sigma=0.2$. On the right is the empirical stationary distribution of the realisation.} \label{fig:fig1}
\end{figure*}

\section{Bifurcation analysis in discrete stochastic systems}\label{sec:bifurc}

We introduce concepts in bifurcation analysis through the normal form of the pitchfork bifurcation, described by the map

\begin{align}
	x(t+1) = rx(t)-x^3(t). \label{eq:bifurc}
\end{align}

Bifurcation analysis concerns the number and stability of fixed points of a system as a function of a parameter -- here $r$. When these change, the system is said to go through a bifurcation, a qualitative change in long term behaviour. Equation \cref{eq:bifurc} has a fixed point at $\xh_0=0$ for all $r$ values. We understand the stability of a fixed point in a map by the magnitude of the derivative (or the eigenvalues of the Jacobian in higher dimensional systems). If it is smaller in absolute value than $1$, the system will move towards the point, otherwise away from it, as can be checked by so-called cobwebbing. As $\frac{\partial rx-x^3}{\partial x} |_{x=\xh_0} = r$, the fixed point at 0 exhibits a change in its stability at $r=1$. Furthermore, for $r>1$ two new fixed points arise, $\xh_{\pm1} = \pm \sqrt{r-1}$. Examining their stability, $\frac{\partial rx-x^3}{\partial x} |_{x=\xh_{\pm1}} = -2r+3$, they are initially stable, becoming unstable at $r>2$. This analysis of the system is summarised in \Cref{fig:fig1}A.

Unfortunately even minuscule amounts of noise destroy this beautiful picture \cite{Crauel1998}. Stochastic bifurcation analysis examines random dynamical systems, and distinguishes between phenomenological and dynamical bifurcations. The former examines profiles of stationary probability densities, whereas the latter examines systems with frozen noise. Numerical analysis is best suited towards finding phenomenological bifurcations, as by definition these cause noticeable changes in system behaviour \cite{Arnold1991}. As experimentally found, the number and stability of the so-called random fixed points often differs from the noiseless system's fixed points at a particular parameter value, thus also changing the point of bifurcation \cite{Diks2006,Kuehn2011,Wang2018}. To illustrate the difficulty of the problem in even our simple example, examine \Cref{fig:fig1} B and C, where a deterministic system already over the point of bifurcation falls back to its former behaviour with added noise.

In known systems, the expected number of fixed points may be given, but in real world measured datasets we have little idea. Therefore we need a model that is capable of automatically determining the number and location of likely fixed points, given little data. In the following section we describe such an approach, based on explicitly parametrising Gaussian Process transition functions in terms of fixed points and corresponding Jacobians, incorporating a form of Automatic Relevance Determination to identify the number of constraints best supported by the data.

\section{Model}\label{sec:model}

In this section, we first introduce the type of datasets we expect for a single value of the bifurcation parameter, which arise frequently in biological studies, and the assumptions that underlie the models used here. We then describe the general machinery of learning a Gaussian Process transition map from data. Finally, in \cref{sec:fp_extension} we describe our extension of this machinery, which enables the fixed points and the corresponding local linearisations to be optimised directly.

Given a dataset of repeated measurements of time-series $\mathsf{Y} = \{y_e^{t,n} \in \mathbb{R}^{D_y}\}_{t \in [0, T]}^{n \in \mathrm{Trials}}$ of $D_y$ measured variables, we aim to model the data using a discrete-time, latent, non-linear, stochastic dynamical system:
\begin{align}
	x_d^{t,n} &= f_d(x_d^{t-1,n}) + \varepsilon_d \\
	y_e^{t,n} &= g_e(x_d^{t,n}) + \eta_e,
\end{align}
where $x_d^{t,n} \in \mathbb{R}^{D_x}$ represents the state of the dynamical system at time $t$ on trial $n$; $f: \mathbb{R}^{D_x} \rightarrow \mathbb{R}^{D_x}$ is the transition function, with $\varepsilon_d \sim \mathcal{N}\left(0,(\sigma^\varepsilon_d)^2\right)$ additive noise; and $g: \mathbb{R}^{D_x} \rightarrow \mathbb{R}^{D_y}$ is the observation function, with $\eta_e \sim \mathcal{N}\left(0,(\sigma^\eta_e)^2\right)$ additive noise.

For simplicity, we only discuss the case where $g$ is a linear function, $g_e(x_d) = C_{ed}x_d$, and $D_x \leq D_y$. We can thus write down the conditional probability
\begin{align}
	 p(y_e^{t,n} \given x_d^{t,n}) = \mathcal{N}\left(y_e^{t,n} - C_{ed}x_d^{t,n},(\sigma^\eta_e)^2\right).
\end{align}

The framework may be extended easily to arbitrary observation functions and likelihoods (including Gaussian Processes models of the observation function), as long as we can compute the probability $p(y_e^{t,n} \given x_d^{t,n})$, estimate its gradient with respect to $x_d^{t,n}$ and infer an approximation to the posterior $p(x_d^{t,n} \given x_d^{t-1,n}, y_e^{t,n})$.

\subsection{Sparse Gaussian process transition function inference and learning}

We represent the scalar-valued stochastic transition functions $f^d$ as independent Gaussian Processes (GP), such that the latent map is $f_d = \left[f^1, f^2, \dots , f^{D_x}\right]^\tr$, and we can easily query the predicted mean and variance of the map for any given input. Unforunately, inferring a full GP from time series observations requires re-estimating the posterior after incorporating every new data point, then reevaluating the influence of previous observations on the current belief state, given the new estimate of the transition functions. This is very much intractable, but choosing a parameterised representation of the function over inducing points de-couples the current estimate of the function from the observations \cite{McHutchon2014} and enables us to evaluate the model.

The process $f_d$ is thus a Sparse Gaussian Process (SGP) parametrised by: a positive definite kernel function $K: \mathbb{R}^{D_x\times N} \otimes \mathbb{R}^{D_x \times M} \rightarrow \mathbb{R}^{N \times M} $; the inducing point locations $z_{dm} \in \mathbb{R}^{D_x \times M}$; and the uncertain values at those locations, represented as random variables drawn from $\mathcal{N}\left(u_{dm},(\sigma^u_m)^21_d\right)$. The ability to represent different noise levels at different locations makes it possible to learn maps with heteroscedastic noise, and may also serve as an Automatic Relevance Determination (ARD)  system to determine the required set of inducing points. Given the parameters, we can evaluate the expected mean and variance of each map at any input value $x_d$:
\begin{align}
	K_{m_1m_2}^{zz} &= K_{m_1m_2}(z_{dm}, z_{dm})+(\sigma^u_m)^2\delta_{mm_1m_2}  \label{eq:inf_start}\\
	\alpha_{dm_1} &= \left[K_{m_1m_2}^{zz}\right]^{-1}u_{dm_2} \\
	\E_f[f^d(x_d)] &= K_{m_1}(x_d, z_{dm})\ \alpha_{dm_1} \\
	\mathrm{Var}_f[f^d(x_d)] &= k(x_d, x_d) -  \nonumber \\ 
	 K_{m_1}(x_d, &z_{dm})\left[K_{m_1m_2}^{zz}\right]^{-1}K_{m_2}(x_d, z_{dm})
\end{align}

We then recover the vector-valued transition function as a collection of the scalar functions, $f_d(x) = f^{d}(x),\quad d=1\dots D_x$. The learning of the transition function thus boils down to the estimation of the inducing point parameters $\theta_{\textrm{ind}} = \{z_{dm},u_{dm},\sigma_m^u\}$. Often the kernel function itself is parametrised, $K = K(\cdot,\cdot \given \theta_{\textrm{k}})$, in which case we may choose to learn the kernel hyperparameters, $\theta_k$, too.

\subsubsection{Inference}

In order to carry out the inference of a latent trajectory, given the data and the current estimate of the transition and observation functions, there are various algorithms available \cite{882463}. We choose to use Assumed Density Filtering \cite{5967674} here, based on empirical performance in tasks similar to ours \cite{McHutchon2014}. 

We represent our belief of the latent state at time $t-1$ as a normal distribution with mean $\mu_d^\toto$ and diagonal covariance matrix $\Sigma_{d_1d_2}^\toto$. We first need to propagate our belief through the transition function, to get an estimate of our updated belief, which we approximate as a normal distribution $\mathcal{N}(\mu_d^\tot, \Sigma_{d_1d_2}^\tot)$, with a non-diagonal covariance. The moments of this distribution may then be computed as:
 \begin{align}
	\begin{split}
	\mu_d^\tot &= \E_{x}[\E_f[f_d(x_d)]] \\
	&= \E_x[K_{m_1}(x_d, z_{dm})]\alpha_{dm_1} \\
	\end{split} \\[10pt]
	\begin{split} \label{eq:inf_end}
		\Sigma_{d_1d_2}^\tot &= \E_{x}[\mathrm{Var}_f[f_d(x_d)]]\delta_{d_1d_2} \\
		&= \E_x[k(x_d, x_d)]\delta_{d_1d_2} - \\ 
		&\quad \big(\E_x[K_{m_1}(x_d, z_{dm})K_{m_2}(x_d, z_{dm})]\  \cdot \\
		&\qquad\qquad\cdot \ \left[K_{m_1m_2}^{zz}\right]^{-1}\big)\delta_{d_1d_2} \\
	\end{split}
\end{align}
where the required expectations may be computed in closed form for linear and Exponentiated Quadratic kernels, and are shown in Appendix B.

Given our belief of the latent state at time $t$, we need to incorporate the data into our belief to obtain the moments. Thanks to our simple linear-Gaussian observation model, and our approximate belief, this can be done exactly:
\begin{align}
	\Sigma^\textrm{fwd}_{e_1e_2} &= C_{e_1d_1} \Sigma_{d_1d_2}^\tot C_{e_2d_2} + (\sigma^\eta_e)^2 \delta_{ee_1e_2} \\[8pt]
	\Sigma^\textrm{back}_{de} &=  \Sigma_{dd_1}^\tot C_{e_1d_1}[\Sigma^\textrm{fwd}_{e_1e}]^{-1}\\[8pt]
	\mu_d^\ttt &= \mu_d^\tot + \Sigma^\textrm{back}_{de_1}\left(y^t_{e_1} - C_{e_1d_1}\mu_{d_1}^\tot \right) \\[8pt]
	\Sigma_{d_1d_2}^\ttt &= \Sigma_{d_1d_2}^\tot - \Sigma^\textrm{back}_{d_1e_1}C_{e_1d_3}\Sigma_{d_3d_2}^\tot
\end{align}

We then approximate the covariance matrix $\Sigma_{d_1d_2}^\ttt$ with its diagonal, as per moment matching, and proceed to carry out the filtering for the next time step, until the complete trajectory has been recovered. The parameters of the latent model we require for inference are $\theta_\textrm{LDS}=\{ C_{ed}, \sigma^\eta_e, \mu^{0,0}_d, \Sigma^{0,0}_{d_1,d_2} \}$.

\subsubsection{Learning}

In order to obtain a good estimate of the parameters of our model, we need to be able to learn them. There are many frameworks available to carry out this estimation \cite{pmlr-v5-titsias09a,Titsias2010,1605.07066}, we chose - again, based on empirical evidence - to use gradient ascent with the exact log marginal likelihood of the above described model as our objective function.
\begin{align}
	\begin{split}
		\mathcal{L}(\theta) &= \log p(Y_{dtn} \given \theta) \\
			&= \sum_{n=1}^N \sum_{t=1}^T \log p(y^{t,n}_d \given y^{1:(t-1),n}_d, \theta) \\
			&= \sum_{n,t} \log p(y^{t,n}_d \given \mu^{t-1,t,n}_d, \Sigma^{t-1,t,n}_{d_1d_2}, \theta) \\
			&= \sum_{n,t}  -\frac{1}{2}\Bigl(D_y\log(2\pi) + \log\left|\Sigma^{\textrm{fwd}, t,n}_{d_1d_2}\right| \\
			&\qquad \qquad+\left(y^{t,n}_{e_1} - C_{e_1d_1}\mu_{d_1}^{t-1,t,n} \right) [\Sigma^{\textrm{fwd}, t,n}]^{-1}_{e_1e_2} \ \cdot \\
			&\qquad\qquad \cdot \ \left(y^t_{e_2} - C_{e_2d_1}\mu_{d_1}^{t-1,t,n} \right)\Bigr)
	\end{split}
\end{align}

Equipped with this objective function, and the fact, that our model is differentiable with respect to all of its parameters, we may optimise our parameters via gradient ascent. We may then iterate inference and learning steps until convergence.

\subsection{Conditioning on fixed points} \label{sec:fp_extension}

We now wish to extend our framework towards a Fixed Point Sparse Gaussian Process, which provides explicit representation of the learned map's fixed points and their local linearisation as parameters of the model.

We may think of a fixed point in a Sparse Gaussian Process as a special inducing point, whose value is tied to its location. Furthermore, to represent and inquire about the stability of said fixed point, we wish to attach derivative observations to that location, representing the local Jacobian. This way we may uncover the location and stability of fixed points in non-linear dynamical systems, via any optimisation algorithm.

The steps we need to go through for the derivation of the system largely follow what has been described in the previous section, with a few extra complications. Let the fixed points be represented as random variables drawn from $\mathcal{N}\left(s_{dp},(\sigma^s_p)^21_d\right)$, where the $\sigma^s_p$ variables may be used by the system to disable unnecessary fixed points, by setting this variance high. Each fixed point is associated with a local $D_x\times D_x$ Jacobian, stored as a tensor $J_{d_1d_2p}$.

\begin{figure*}[t] 
\centering
\includegraphics[width=0.8\linewidth]{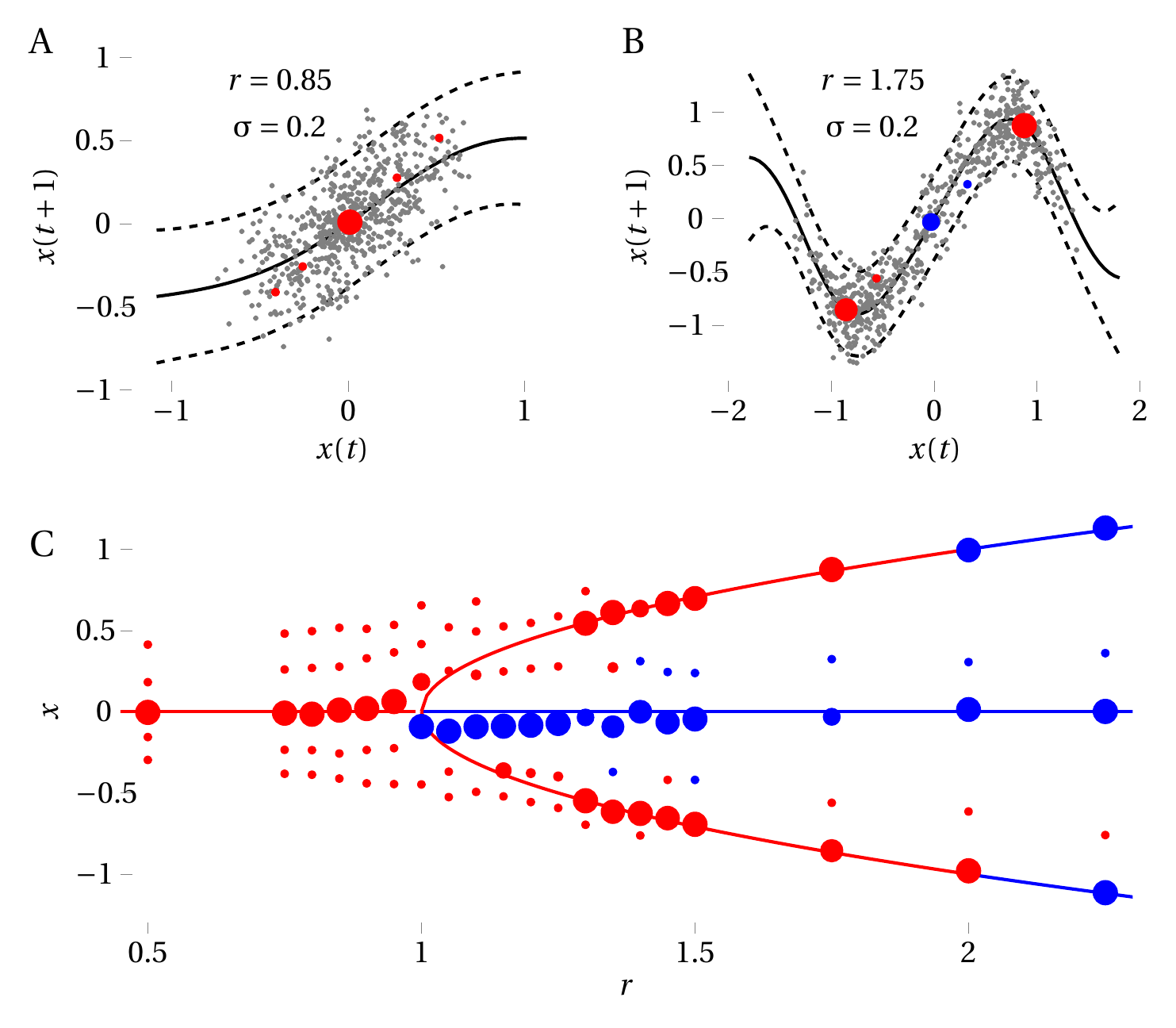}
\caption{Stochastic bifurcation experiment. {\it A.} The posterior fit of the Fixed Point Sparse Gaussian Process on simulated data from \cref{eq:bifurc_stochastic}. Small grey points represent the training data. The solid line is the posterior mean, with the dashed lines indicating the posterior standard deviation at $2\sigma_\textrm{post}$. The red and blue circles are the learned fixed point parameter locations, with colour indicating stability and size indicating the strength of belief in the fixed point. {\it B.} Same as {\it A}, with $r=1.75$, after the bifurcation. {\it C.} Bifurcation plot of the stochastic system. The solid lines represent the noiseless bifurcation, and are equivalent to \Cref{fig:fig1}A. The circles represent the learned fixed point locations [y axis position], for a dataset simulated at a particular value of $r$ [x axis position]. Marker colour and size as described in {\it A}.}
\label{fig:fig2}
\end{figure*}

\begin{figure*}[t]
\centering
\includegraphics[width=0.7\linewidth]{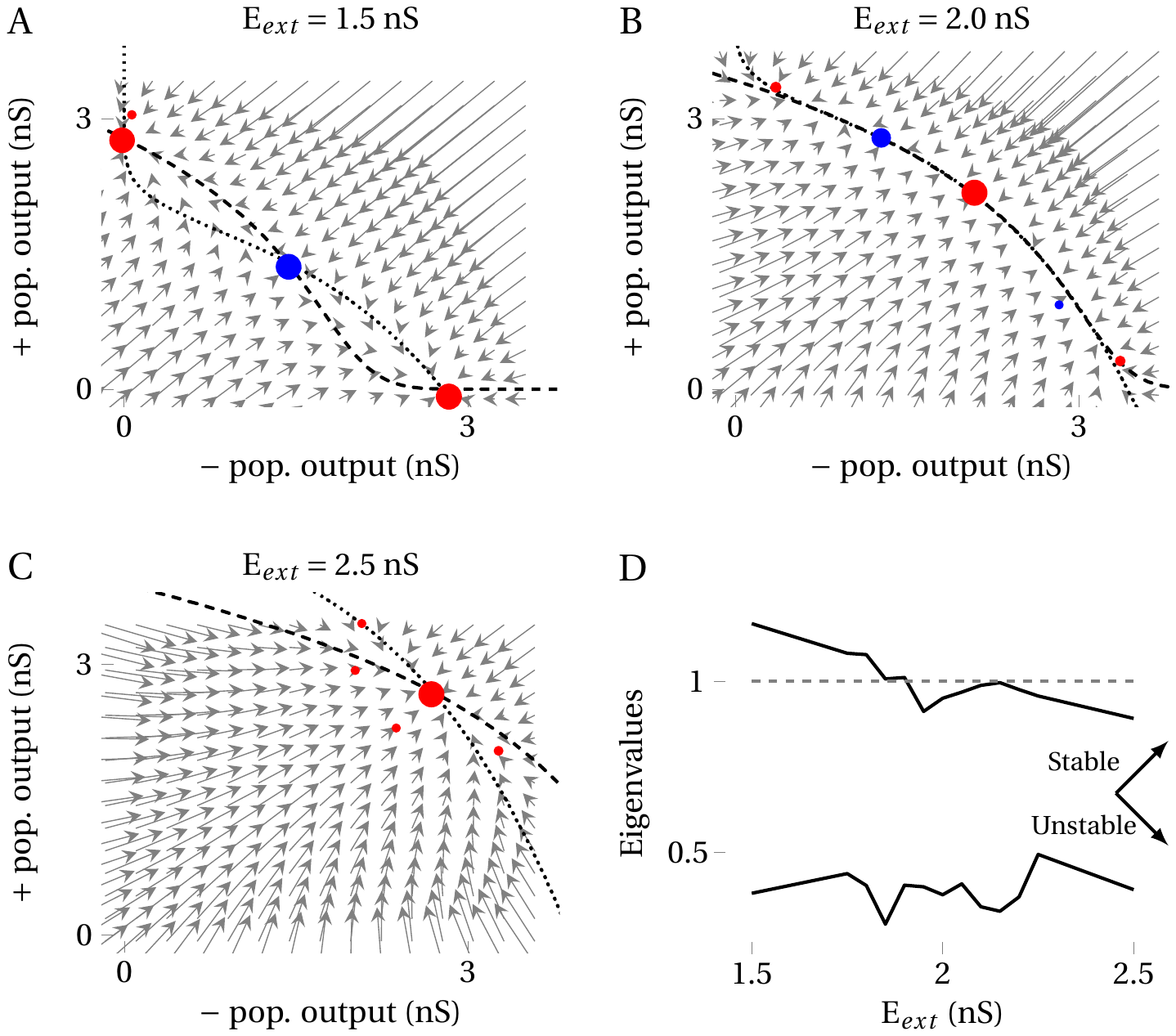}
\caption{Mutual inhibition experiment. {\it A-C.} The posterior fit of the Fixed Point Sparse Gaussian Process on simulated data from \Crefrange{eq:pos}{eq:neg} for various $E_{ext}$ excitatory input values. The axis are the population outputs, and the learned posterior mean is described by the grey arrow. The dashed and dotted lines are the true nullclines of the system, for the positive and negative populations, respectively. The circles represent the learned fixed point locations, with red and blue colours indicating stable and unstable locations. The size of the marker corresponds to strength of belief in the fixed point. {\it D.}  The learned eigenvalues of the central fixed point, with varying $E_{ext}$. The inset indicates the average eigenvector direction over all experiments, with the more unstable direction pointing along the nullclines.
} \label{fig:fig3}
\end{figure*}

We can thus extend the parameter set describing our current belief of the map, $\theta_{\textrm{ind+fp}} = \{z_{dm},u_{dm},\sigma_m^u, s_{dp}, \sigma^s_p, J_{d_1d_2p} \}$. We use this set of parameters to carry out the inference, requiring us to revisit \crefrange{eq:inf_start}{eq:inf_end}. The main change comes from the fact, that we wish to use the derivative observations attached to our fixed points during inference. This requires establishing a derivative Gaussian Process, whose kernel function is given by the derivative of the original kernel function with respect to both arguments
\begin{align}
	 K_{nm,d^1d^2}& = \nabla_1\nabla_2 K_{nm}(s_{d^1n}, s_{d^2m} \,|\, \theta_k),
\end{align}
resulting in a fourth order tensor. It is useful to define a block-structure matrix version $K_{nm,d^1d^2}  \triangleq \bK_{(nd^1)(md^2)}$. We may similarly define the cross-covariance between a normal and a derivative process, as taking the derivative with only to the respective argument of the original kernel function:
\begin{align}
	 K_{nm,d^1}& = \nabla_1 K_{nm}(s_{d^1n}, x_{d^2m} \,|\, \theta_k) \\
 	 K_{nm,d^2}& = \nabla_2 K_{nm}(x_{d^1n}, s_{d^2m} \,|\, \theta_k)
\end{align}

Equipped with these processes, we may re-write the predictive \crefrange{eq:inf_start}{eq:inf_end} as the block matrices
\[
K^\textrm{block}_{ab}=
\left[
\begin{array}{c c c}
\vspace{0.5cm}
K^{zz}_{m_1m_2} & K^{zs}_{mp} & K^{zJ}_{m(pd)} \\
\vspace{0.5cm}
K^{sz}_{m_1m_2} & K^{ss}_{p_1p_2} & K^{sJ}_{p_1(p_2d)} \\

K^{Jz}_{(pd)m} & K^{Js}_{(p_1d)p_2} & K^{JJ}_{(p_1d_1)(p_2d_2)}
\end{array}
\right]
\]

\[
\alpha_{db_1} = \left[K^\textrm{block}_{b_1b_2}+
\left[
\begin{array}{c}
\sigma^u_{m} \\
\sigma^s_{p}\\
0_{(pd)}
\end{array}
\right]\delta_{bb_1b_2}
\right]^{-1}
\left[
\begin{array}{c}
u_{dm} \\
s_{dp}\\
J_{d(pd_2)}
\end{array}
\right]
\]

\[
K_b^\textrm{pred} = 
\left[
\begin{array}{c c c}
K^{xz}_{m} & K^{xs}_{p} & K^{xJ}_{(pd)}
\end{array}
\right]
\]

Resulting in the predictive moments for a noiseless input:
\begin{align}
\E_f[f^d(x_d)] &= K_b^\textrm{pred}\alpha_{db} \\
\mathrm{Var}_f[f^d(x_d)] &= k(x_d, x_d) - \nonumber \\
&\quad K_a^\textrm{pred} \left[K^\textrm{block}_{ab}\right]^{-1} K_b^\textrm{pred},
\end{align}

where $a$ and $b$ represent the concatenation of indices $[m, p, (pd)]$. For the inference we still need to consider propagating beliefs represented as Gaussian random variables, thus requiring to compute $\E_{x}[\E_f[f_d(x_d)]]$ and $E_{x}[\mathrm{Var}_f[f_d(x_d)]]$. For the linear and the Exponentiated Quadratic kernels these expressions are available in closed form, although care must be taken during the computations, as for the latter term, we need to take the expectation of kernel products between regular and derivative processes. In the general case these expectation are to be approximated numerically \cite{Girard2004}. For the detailed computations, see Appendix B.

The learning does not change significantly, our objective function remains the same, and our operations remain differentiable with respect to all parameters.

\section{Experiments}\label{sec:experiments}

Equipped with fully described model, we are ready to test it. As our first example we are going to return to the well-studied example described in \cref{sec:bifurc}. Finally, we study changes of fixed point pattens in an influential model of mutually inhibiting neural populations during decision making \cite{Machens2005}.

\subsection{Stochastic pitchfork bifurcation}

We may now write down the stochastic version of \cref{eq:bifurc}:
\begin{align}
	x(t+1) = rx(t)-x^3(t) + \epsilon_t \label{eq:bifurc_stochastic}
\end{align}

where $\epsilon_t \sim \mathcal{N}\left(0, \sigma^2\right)$ iid. We examine, how varying $r$ affects the learned fixed points. We trained the system using 32 trials, lasting 20 time steps each, with $\sigma=0.2$, and the initial condition $x(0) \sim \mathcal{N}(0,0.001^2)$. Note that this is less data in total than shown in \Cref{fig:fig1} C. We then fit our model to the data with 16 inducing points and the overestimated 5 fixed points, letting the ARD formulation determine the number of fixed points present in the system.

We first confirmed, that the method indeed captures the available data very well for various values of the bifurcation parameter $r$, as shown in, \Cref{fig:fig2} A and B. 

We then create the bifurcation plot, \Cref{fig:fig2}C, based on the learned parameter values. The fixed points identified truthfully track the expected location and stability, as well as successfully recovering the true number of fixed points. Consistently with previous finding on similar systems \cite{Diks2006}, we indeed find that noise shifts the bifurcation towards larger values of $r$, and when the distances of the noiseless fixed points are on the order of the noise, the random fixed points are not detectable from data.

\subsection{Mutually inhibiting neural populations}

Having recovered previous results with our highly flexible system, we now turn our attention to a system closer to the data-analysis-in-the-wild type problems we aimed to solve.

The data used comes from a simulated system, but one that was optimised to match the behaviour of measured neural population. For more details about the experiment and the simulation, read the excellent paper from Machens, Romo and Brody \yrcite{Machens2005}. 

\begin{figure}[H]
\begin{center}
\begin{tikzpicture}[every plot/.append style={very thick}]
	\node[circle, draw=black, fill=none] at (4,2) (plus){\large $+$};
	\node[circle, draw=black, fill=none] at (4,0) (minus){\large $-$};
	\draw	(minus)to[out=60,in=-60] (plus);
	\draw	(minus)to[out=130,in=-130] (plus);
	\filldraw (4.23,1.63) circle (3pt);
	\filldraw (3.77,0.37) circle (3pt);
	\node[draw=none, fill=none] at (2.6,1) (E){\large E$_{ext}$};

	\draw	(E)to[out=0,in=180] (minus);
	\draw	(E)to[out=0,in=180] (plus);

\end{tikzpicture}
\caption{Circuit diagram of externally stimulated mutually inhibiting neural populations}
\end{center}
\end{figure}
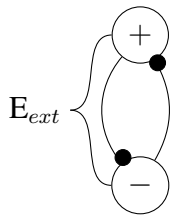

In short the system consists of an external excitation and so-called negative and positive populations. In the current study we do not take into account possible differential inputs to the populations, so for our purposes the model is completely symmetric. We slightly reformulate the system equations to match the language used throughout the paper. The simulated system is thus governed by the map

\begin{align}
Pos(t+1) &= Pos(t) + f(-\omega_I Neg(t) + E_{ext}) + \epsilon^{t,0} \label{eq:pos}\\
Neg(t+1) &= Neg(t) + f(-\omega_I Pos(t) + E_{ext}) + \epsilon^{t,1},\label{eq:neg}
\end{align}

where $\epsilon^t \sim \mathcal{N}\left(0, \sigma^2\right)$ is iid noise, $\omega_I$ is the inhibitory weight of the populations, and importantly, $f$ is a numerically optimised function, which defines the nullclines of the system and gives rise to interesting system behaviour, as the external driving input $E_{ext}$ - our bifurcation parameter - is varied.

The nullclines generally cross one another, giving rise to fixed points, whose stability depend on the angle of crossing. In particular, the original study finds that in one extreme, two stable fixed points are created far from one another, and may be used for decision making, whereas during the other, there is only a single fixed point, useful for loading new information robustly. The special case is, that for a special value of the external input, the nullclines were designed such that they completely overlap, creating a line attractor, which is used for maintaining system state.

Note that this model was carefully designed to behave so, and our goal here is to estimate such behaviour purely from data generated from the model. We simulated 60 two-dimensional trajectories, each consisting of 80 time steps, corresponding to 2 seconds of experimental data, sampled at 25 msec steps, and injecting additive noise at every time step, with $\sigma=0.1$.

Examining the results in \Cref{fig:fig3} A-C the system behaviour was very well captured by the estimated fixed points, including the number, location and stability of the points. Although the stability of the system is inferred correctly, if we examine the eigenspectrum of the central point in \Cref{fig:fig3}D, we can indeed follow its stabilisation from a saddle with one stable and one unstable direction to a mostly stable fixed point.

\section{Discussion} \label{sec:discussion}

Studying real systems, especially in biological experiments, where we have little knowledge of the governing equation is a hard but ubiquitous problem. In the current study we designed an algorithm aimed at the study of random dynamical systems measured at discrete time, in which we can modify or measure a variable influencing the system. We provide a simple analysis of a very complicated question, namely when does the system behaviour change qualitatively as a system parameter varies. This question broadly encompasses large fields of research, particularly cell biology and neuroscience.

All analyses comes with limitations of course. The current work is aimed at fixed point bifurcations, at the moment we can not sufficiently describe limit cycles or more than zero-dimensional attractors, beyond a rudimentary approximation in the form of aligned fixed points. Another issue that often comes up in learning dynamical systems is that of global stability. Fortunately Gaussian Processes with many types of kernels functions decay rapidly away from data, so likely our learned map would push us towards the origin as soon as a new input wanders away from observations.

Our core contribution is the Fixed Point Sparse Gaussian Process formulation, in which fixed points appear explicitly as parameters of the model fit to data, and may thus be identified directly by parameter optimisation methods. This core idea has many potential extensions highlighted throughout the paper, including extensions to non-linear maps from latent dynamics to observations, or to hierarchical (deep) Gaussian Processes models of the transition map itself, which are capable of capturing more complex and less smooth structures.

Furthermore, combining this powerful stochastic representation with the ability to robustly identify fixed points in unknown systems may indeed bring further effort into stochastic bifurcation analysis, an exciting and very powerful methodology, still in its infancy.

\bibliographystyle{icml2018}
\bibliography{var_sgp_calcium}

\begin{thebibliography}{17}
\providecommand{\natexlab}[1]{#1}
\providecommand{\url}[1]{\texttt{#1}}
\expandafter\ifx\csname urlstyle\endcsname\relax
  \providecommand{\doi}[1]{doi: #1}\else
  \providecommand{\doi}{doi: \begingroup \urlstyle{rm}\Url}\fi

\bibitem[Arnold \& Crauel(1991)Arnold and Crauel]{Arnold1991}
Arnold, Ludwig and Crauel, Hans.
\newblock Random dynamical systems.
\newblock pp.\  1--22, 01 1991.

\bibitem[Bui et~al.(2016)Bui, Yan, and Turner]{1605.07066}
Bui, Thang~D., Yan, Josiah, and Turner, Richard~E.
\newblock A unifying framework for gaussian process pseudo-point approximations
  using power expectation propagation, 2016.

\bibitem[Crauel \& Flandoli(1998)Crauel and Flandoli]{Crauel1998}
Crauel, Hans and Flandoli, Franco.
\newblock {Additive noise destroys a pitchfork bifurcation}.
\newblock \emph{J. Dyn. Differ. Equations}, 10\penalty0 (2):\penalty0 259--274,
  1998.
\newblock ISSN 1040-7294.
\newblock \doi{10.1023/A:1022665916629}.
\newblock URL \url{http://link.springer.com/article/10.1023/A:1022665916629}.

\bibitem[Diks(2006)]{Diks2006}
Diks, C.
\newblock {A weak bifurcation theory for discrete time stochastic dynamical
  systems}.
\newblock \emph{Tinbergen Inst. Discuss. Pap. No. 06-}, \penalty0
  (June):\penalty0 1--30, 2006.
\newblock URL
  \url{http://papers.ssrn.com/sol3/papers.cfm?abstract{\_}id=901422}.

\bibitem[Girard(2004)]{Girard2004}
Girard, Agathe.
\newblock {Approximate methods for propagation of uncertainty with Gaussian
  process models}.
\newblock \emph{Ph.D. Thesis}, \penalty0 (October), 2004.
\newblock URL
  \url{http://citeseerx.ist.psu.edu/viewdoc/download?doi=10.1.1.66.1826{\&}rep=rep1{\&}type=pdf}.

\bibitem[Kuehn(2011)]{Kuehn2011}
Kuehn, Christian.
\newblock {A mathematical framework for critical transitions: Bifurcations,
  fastslow systems and stochastic dynamics}.
\newblock \emph{Phys. D Nonlinear Phenom.}, 240\penalty0 (12):\penalty0
  1020--1035, 2011.
\newblock ISSN 01672789.
\newblock \doi{10.1016/j.physd.2011.02.012}.

\bibitem[Machens et~al.(2005)Machens, Romo, and Brody]{Machens2005}
Machens, Christian~K., Romo, Ranulfo, and Brody, Carlos~D.
\newblock {Flexible control of mutual inhibition: A neural model of
  two-interval discrimination}.
\newblock \emph{Science (80-. ).}, 307\penalty0 (5712):\penalty0 1121--1124,
  2005.
\newblock ISSN 00368075.
\newblock \doi{10.1126/science.1104171}.
\newblock URL \url{http://www.sciencemag.org/cgi/doi/10.1126/science.1104171}.

\bibitem[McHutchon(2014)]{McHutchon2014}
McHutchon, Andrew.
\newblock {Nonlinear Modelling and Control using Gaussian Processes}.
\newblock \penalty0 (August), 2014.
\newblock URL \url{.}

\bibitem[Nonnenmacher et~al.(2017)Nonnenmacher, Turaga, and
  Macke]{Nonnenmacher2017}
Nonnenmacher, Marcel, Turaga, Srinivas~C, and Macke, Jakob~H.
\newblock {Extracting low-dimensional dynamics from multiple large-scale neural
  population recordings by learning to predict correlations}.
\newblock \emph{Adv. Neural Inf. Process. Syst.}, \penalty0 (Nips), 2017.

\bibitem[Park et~al.(2015)Park, Bohner, and Macke]{Park2015}
Park, Mijung, Bohner, Gergo, and Macke, Jakob~H.
\newblock {Unlocking neural population non-stationarities using hierarchical
  dynamics models}.
\newblock \emph{Adv. Neural Inf. Process. Syst.}, pp.\  145--153, 2015.

\bibitem[Ramakrishnan et~al.(2011)Ramakrishnan, Ertin, and Moses]{5967674}
Ramakrishnan, N., Ertin, E., and Moses, R.~L.
\newblock Assumed density filtering for learning gaussian process models.
\newblock In \emph{2011 IEEE Statistical Signal Processing Workshop (SSP)},
  pp.\  257--260, June 2011.
\newblock \doi{10.1109/SSP.2011.5967674}.

\bibitem[Sussillo \& Barak(2013)Sussillo and Barak]{Sussillo2013}
Sussillo, David and Barak, Omri.
\newblock {Opening the Black Box: Low-Dimensional Dynamics in High-Dimensional
  Recurrent Neural Networks}.
\newblock \emph{Neural Comput.}, 25\penalty0 (3):\penalty0 626--649, 2013.
\newblock ISSN 0899-7667.
\newblock \doi{10.1162/NECO_a_00409}.
\newblock URL
  \url{http://www.mitpressjournals.org/doi/10.1162/NECO{\_}a{\_}00409}.

\bibitem[Sussillo et~al.(2016)Sussillo, Jozefowicz, Abbott, and
  Pandarinath]{Sussillo2016}
Sussillo, David, Jozefowicz, Rafal, Abbott, L.~F., and Pandarinath, Chethan.
\newblock {LFADS - Latent Factor Analysis via Dynamical Systems}.
\newblock \emph{arXiv}, 2016.
\newblock URL \url{http://arxiv.org/abs/1608.06315}.

\bibitem[Titsias(2009)]{pmlr-v5-titsias09a}
Titsias, Michalis.
\newblock Variational learning of inducing variables in sparse gaussian
  processes.
\newblock In van Dyk, David and Welling, Max (eds.), \emph{Proceedings of the
  Twelth International Conference on Artificial Intelligence and Statistics},
  volume~5 of \emph{Proceedings of Machine Learning Research}, pp.\  567--574,
  Hilton Clearwater Beach Resort, Clearwater Beach, Florida USA, 16--18 Apr
  2009. PMLR.
\newblock URL \url{http://proceedings.mlr.press/v5/titsias09a.html}.

\bibitem[Titsias \& Lawrence(2010)Titsias and Lawrence]{Titsias2010}
Titsias, Michalis and Lawrence, Neil.
\newblock {Bayesian Gaussian Process Latent Variable Model}.
\newblock \emph{Artif. Intell.}, 9:\penalty0 844--851, 2010.
\newblock ISSN 0899-7667.
\newblock \doi{10.1162/089976699300016331}.
\newblock URL \url{http://eprints.pascal-network.org/archive/00006343/}.

\bibitem[Wan \& Merwe(2000)Wan and Merwe]{882463}
Wan, E.~A. and Merwe, R. Van~Der.
\newblock The unscented kalman filter for nonlinear estimation.
\newblock In \emph{Proceedings of the IEEE 2000 Adaptive Systems for Signal
  Processing, Communications, and Control Symposium (Cat. No.00EX373)}, pp.\
  153--158, 2000.
\newblock \doi{10.1109/ASSPCC.2000.882463}.

\bibitem[Wang et~al.(2018)Wang, Chen, and Duan]{Wang2018}
Wang, Hui, Chen, Xiaoli, and Duan, Jinqiao.
\newblock {A Stochastic Pitchfork Bifurcation in Most Probable Phase
  Portraits}.
\newblock pp.\  1--9, 2018.
\newblock URL \url{http://arxiv.org/abs/1801.03739}.

\end{thebibliography}

\end{document}